\documentclass[preprint,12pt]{elsarticle}
\usepackage{graphicx} 
\usepackage{amssymb}
\usepackage{amsthm}
\usepackage{graphics}
\usepackage{epsfig}
\usepackage{amssymb}
\usepackage{hyperref}
\usepackage{graphicx}
\usepackage{subfigure}
\usepackage{float}
\usepackage{color}
\usepackage{tikz-network}
\usepackage{tikz}
\usepackage{tikz,pgfplots}
\usetikzlibrary{patterns}
\usepackage[ruled,vlined]{algorithm2e}
\usepackage{caption}
\DeclareCaptionLabelFormat{cont}{#1~#2\alph{ContinuedFloat}}
\usepackage{fullpage}
\usepackage{amsmath}
\usepackage{amsfonts}
\usepackage{enumerate}
\usepackage{enumitem} 
\usepackage{tikz}
\usepackage{array}
\usepackage{bigdelim}
\usepackage{nicematrix}
\usepackage{mathtools}
\usepackage{blkarray, bigstrut}
\usepackage{gauss}
\usepackage{xcolor}
\usepackage{tabularray}
\UseTblrLibrary{booktabs}

\newcommand{\pf}{\noindent {\bf Proof: }}

\newtheorem*{theoremaux}{Theorem \theoremauxnum}
\gdef\theoremauxnum{1}

\newtheorem{lemma}{\bf Lemma}[section]
\newtheorem{theorem}{\bf Theorem}[section]

\newtheorem{proposition}[lemma]{\bf Proposition}
\newtheorem{corollary}[lemma]{\bf Corollary}
\newtheorem{definition}{\bf Definition}[section]
\theoremstyle{definition}
\newtheorem{remark}{\bf Remark}[section]
\newtheorem{example}{\bf Example}[section]

\journal{~}

\begin{document}

\begin{frontmatter}
\title{{Algebraic degree of Cayley colour graphs}}





 \author{Sauvik Poddar}
 \ead{sauvikpoddar1997@gmail.com}



%
 \address{Department of Mathematics, Presidency University, 86/1 College Street, Kolkata 700073, India}

\begin{abstract}
The splitting field of a graph $\Gamma$ with respect to a square matrix $M$ associated with $\Gamma$, is the smallest field extension over the field of rationals $\mathbb{Q}$ that contains all the eigenvalues of $M$. The degree of the extension is called the algebraic degree of $\Gamma$ with respect to $M$. In this paper, we completely determine the splitting field of the adjacency matrix of the Cayley colour graph $\operatorname{Cay}(G,f)$ on a finite group $G$, associated with a class function $f:G\to\mathbb{Q}$ and compute its algebraic degree, which generalize the main results of Wu et al. Moreover, we study the relation between the algebraic integrality of two Cayley colour graphs, and deduce the fact that the algebraic degree and distance algebraic degree of a normal Cayley graph are same, generalizing a result of Zhang et al.
\end{abstract}

\end{frontmatter}

\section{Introduction}

Let $\Gamma=(V,E)$ be a graph with vertex set $V$ and edge set $E\subseteq V\times V$.
The \textit{adjacency polynomial} of $\Gamma$ is the characteristic polynomial of its adjacency matrix $A(\Gamma)$, whose zeros are said to be the \textit{adjacency eigenvalues} (or simply, \textit{eigenvalues}) of $\Gamma$. A graph $\Gamma$ is said to be \textit{integral} if all its adjacency eigenvalues are integers. The quest of characterizing integral graphs was first initiated by Harary and Schwenk \cite{harary1974graphs} in $1973$. As a generalization, in $2018$, M{\"o}nius et al. \cite{monius2018graphs} introduced the notion of algebraic degree.  

\begin{definition}\cite{monius2018graphs}
Given a graph $\Gamma$, the algebraic degree, $\operatorname{Deg}(\Gamma)$, is the dimension of the
splitting field of the adjacency polynomial of $\Gamma$ over the field $\mathbb{Q}$. 
\end{definition}

\begin{definition}\cite{abdollahi20242}
A graph $\Gamma$ is called $d$-integral if $\operatorname{Deg}(\Gamma)=d$. 
\end{definition}

Clearly, a graph $\Gamma$ is integral if and only if $\operatorname{Deg}(\Gamma)=1$. 
Following the notion of algebraic degree of the adjacency matrix of a graph, one can generalize this idea for any square matrix associated with a graph. Let $\Gamma$ be a graph and let $M(\Gamma)$ be some square matrix associated with $\Gamma$. We call the eigenvalues of $M(\Gamma)$ as the \textit{$M$-eigenvalues} of $\Gamma$. The spectrum of $M(\Gamma)$ is called the \textit{$M$-spectrum} of $\Gamma$, which is the multiset of its $M$-eigenvalues and is denoted by
$$\operatorname{Spec}(M(\Gamma))=\begin{pmatrix}
\lambda_1 & \lambda_2 & \cdots & \lambda_k\\
m_1 & m_2 & \cdots & m_k
\end{pmatrix},$$
where the first row denotes distinct $M$-eigenvalues of $\Gamma$ and the second row denotes their respective multiplicities. A graph $\Gamma$ is said to be \textit{$M$-rational} (resp. \textit{$M$-integral}) if all of its $M$-eigenvalues are rational numbers (resp. integers). Let $K$ be a number field. A graph $\Gamma$ is said to be \textit{$M$-algebraically integral} over $K$ if all of its $M$-eigenvalues are in $K$. The \textit{$M$-splitting field} of $\Gamma$, denoted by $\mathbb{SF}_M(\Gamma)$ is the smallest field extension over $\mathbb{Q}$ which contains all the $M$-eigenvalues of $\Gamma$. The extension degree $[\mathbb{SF}_M(\Gamma):\mathbb{Q}]$ is called the \textit{$M$-algebraic degree} of $\Gamma$ and is denoted by $\operatorname{Deg}_M(\Gamma)$. A graph $\Gamma$ is said to be \textit{$d$-$M$-integral} if $\operatorname{Deg}_M(\Gamma)=d$. In particular, when we consider the eigenvalues, spectrum, splitting field and algebraic degree of $\Gamma$, we mean with respect to the adjacency matrix. For distance matrix $D(\Gamma)$, we use the terms distance eigenvalues, distance spectrum, distance splitting field and distance algebraic degree.

Let $G$ be a finite group with identity $1$ and $S\subseteq G\setminus\lbrace{1}\rbrace$ be an inverse-closed subset, i.e., $S=S^{-1}$, where $S^{-1}\coloneqq\lbrace{s^{-1}~|~s\in S}\rbrace$. A \textit{Cayley graph} $\operatorname{Cay}(G,S)$ of G with respect to $S$ is a graph with $G$ as the set of vertices and two distinct vertices $u,v\in G$ are adjacent in $\operatorname{Cay}(G,S)$ if and only if $uv^{-1}\in S$. The set $S$ is said to be the \textit{connection set} for
$\operatorname{Cay}(G,S)$ and the size of $S$ is called the \textit{valency} of the Cayley graph. A set $S\subseteq G$ is said to be \textit{normal} in $G$ if $gSg^{-1}=S$ for all $g\in G$. A Cayley graph is said to be \textit{normal} if the connection set is normal. Cayley graphs on cyclic groups $\mathbb{Z}_n$ are called \textit{circulant graphs}.

We now consider the multigraph version. A \textit{multiset} is a set $S$ together with multiplicity function $m_S:S\to\mathbb{N}$, where $m_S(x)$ is a positive integer for every $x\in S$, counting the number of times that $x$ occurs in the multiset. We set $m_S(x)=0$ for $x\notin S$. Let $S\subseteq G\setminus\lbrace{1}\rbrace$ be a multiset. $S$ is said to be \textit{inverse-closed} if $m_S(x)=m_S(x^{-1})$ for every $x\in S$. The \textit{Cayley multigraph} on $G$ with respect to the inverse-closed multiset $S$, denoted by $\operatorname{Cay}(G,S)$, is defined to be the multigraph
with vertex set $G$ such that the number of edges joining $x,y\in G$ is equal to $m_S(xy^{-1})$. The adjacency
matrix of $\operatorname{Cay}(G,S)$ is the matrix whose $(x,y)$-entry is $m_S(xy^{-1})$. A Cayley multigraph is said to be \textit{normal} if the connection multiset is normal. In the special case when $m_S(x)=1$ for every $x\in S$, $\operatorname{Cay}(G,S)$ is a Cayley
graph in the usual sense. 

In recent years, much of the works in spectral graph theory is devoted to characterizing integral and distance integral Cayley graphs \cite{alperin2012integral,cheng2019integral,cheng2023integral,guo2019integral,
huang2021distance,huang2021integral,lu2018integral,so2006integral}. Bridges and Mena \cite{bridges1982rational} characterized a necessary and sufficient condition for the integrality of Cayley multigraphs over abelian groups, while DeVos et. al \cite{devos2013integral} addressed the same problem for Cayley multigraphs over Hamiltonian groups. As a natural generalization of integral graphs, the algebraically integral graphs have drawn plenty of attention in recent times \cite{godsil2025integral,li2013circulant,li2020method}.
Beyond this, some further generalizations, such as the splitting fields and algebraic degrees of some family of Cayley graphs, have been studied of late. For instance, M{\"o}nius has studied the algebraic degrees of circulant graphs of prime order \cite{monius2020algebraic} and later generalized those results by determining the splitting fields and the algebraic degrees of circulant graphs of any order $n$ \cite{monius2022splitting}. Lu and M{\"o}nius \cite{lu2023algebraic} determined the splitting fields and algebraic degrees of Cayley graphs over abelian groups and dihedral groups. Wu et al. \cite{wu2024splitting} have determined the splitting field, distance splitting field, algebraic degree and distance algebraic degree of normal Cayley graphs. Taking a similar approach to \cite{abdollahi20242}, where the authors have characterized the $2$-integral Cayley graphs over abelian groups with valency $2,3,4$ and $5$, Zhang et al. \cite{zhang20262} have provided the characterization for $2$-distance integral Cayley graphs. Some other studies also include computation of algebraic degrees of Cayley hypergraphs \cite{sripaisan2022algebraic}, quasi-abelian semi-Cayley digraphs \cite{wang2024algebraic} and $n$-Cayley digraphs \cite{li2025algebraic}. To get a more compact overview and further background information about the eigenvalues and integrality of Cayley graphs, one may refer to the survey article by Liu and Zhou \cite{liu2022eigenvalues}.


\subsection{\textbf{Our contribution}}

\medskip

This paper mainly focuses on Cayley colour graphs with an associated function $f:G\to\mathbb{Q}$ (introduced in Section \ref{prelim-section}), which serves as a generalization of Cayley multigraphs. In Section \ref{sp-field-alg-deg-Cay-col-graph-section}, we completely determine the splitting field and algebraic degree of Cayley colour graphs where the underlying function $f$ is a class function, which generalizes the main results of Wu et al. [\cite{wu2024splitting}, Theorem $3.1$ and Theorem $4.1$]. In section \ref{alg-deg-dist-alg-deg-normal-Cay-graph-section}, we derive the splitting field and algebraic degree of normal Cayley multigraphs and obtain the distance splitting field and distance algebraic degree of normal Cayley graphs, as a corollary to our main theorem (Theorem \ref{splitting-field-alg-deg-Cay-col-graph-thm}). Finally in Section \ref{alg-int-Cay-col-graph-section}, we study the relation between the algebraic integrality of Cayley colour graphs and show that a normal Cayley graph is $d$-integral if and only if it is $d$-distance integral, generalizing a result of Zhang et al. \cite{zhang20262}.

\section{Preliminaries}\label{prelim-section}

We recall some basics of representation and character theory of a finite group. For more information about representation theory and character theory of finite groups, one can refer to \cite{isaacs1994character,serre1977linear,steinberg2012representation}.

Let $G$ be a finite group and $V$ be a finite dimensional vector space over $\mathbb{C}$. A \textit{representation} $(\rho,V)$ of $G$ is a group homomorphism $\rho:G\to GL(V)$, where $GL(V)$ denotes the group of all invertible linear transformations of $V$. We can simply denote a representation as $\rho$ if $V$ is already understood from the context. The \textit{degree} of $\rho$, denoted by $\deg(\rho)$ is $\dim(V)$, the dimension of $V$. 



The character $\chi_\rho:G\to\mathbb{C}$ of $\rho$ is a map defined as $\chi_{\rho}(g)=\operatorname{Tr}(\rho(g))$,
where $\operatorname{Tr}(\rho(g))$ is the trace of the representation matrix of $\rho(g)$ with respect to some basis of $V$. The \textit{degree} of the character $\chi_\rho$, $\deg(\chi_{\rho})$ is the degree of $\rho$, and is equal to $\chi_\rho(1)$. A subspace $W\le V$ is said to be \textit{$G$-invariant} if $\rho(g)w\in W$ for each $g\in G$ and $w\in W$. 
Obviously, $\{0\}$ and $V$ are always $G$-invariant subspaces, which are called trivial. We say that $\rho$ is an \textit{irreducible representation} and $\chi_\rho$ an \textit{irreducible character} of $G$, if $V$ has no non-trivial $G$-invariant subspaces.
We denote by $\operatorname{Irr}(G)$, the complete set of inequivalent irreducible (complex) characters of $G$. 
A function $f:G\to\mathbb{C}$ is called a \textit{class function} if $f(ghg^{-1})=f(h)$, for all $g,h\in G$ or equivalently, if $f$ is constant on each conjugacy class of $G$. The characters of a finite group $G$ are class functions and the set of all inequivalent, irreducible characters of $G$ forms an orthogonal basis for $Z(L(G))$, where $Z(L(G))$ is the centre of the group algebra $L(G)=\mathbb{C}^G=\lbrace{f\mid f:G\to\mathbb{C}}\rbrace$.

\begin{lemma}(\cite{steinberg2012representation}, Corollary $4.1.10$)\label{A-diagonalizable-lemma}
Let $A\in GL_m(\mathbb{C})$ be a matrix of finite order. Then $A$ is
diagonalizable. Moreover, if $A^n=I$, then the eigenvalues of $A$ are $n$th-roots of unity.
\end{lemma}

Let $\zeta_n=e^{\frac{2\pi i}{n}}$ denote the primitive $n$-th root of unity. The following is a well-known result in character theory. However, we give a proof for the sake of completeness.

\begin{lemma}\label{chi-alg-int-lemma}
Let $G$ be a group of order $n$ and $\chi$ be a character of $G$. Then $\chi(g)\in\mathbb{Q}(\zeta_n)$.
\end{lemma}
\pf Let $\rho:G\to GL_m(\mathbb{C})$ be a representation affording the character $\chi$. As $g^n=1$, $(\rho(g))^n=I$. By Lemma \ref{A-diagonalizable-lemma}, $\rho(g)$ is diagonalizable with eigenvalues $\gamma_1,\ldots,\gamma_m$, which are $n$th-roots of unity. Since $\gamma_i\in\mathbb{Q}(\zeta_n)$ for all $i\in\lbrace{1,\ldots,m}\rbrace$ and $\chi(g)=\operatorname{Tr}(\rho(g))=\gamma_1+\cdots+\gamma_m$, the result follows.\qed

\medskip

For a finite group $G$ and a function $f:G\to\mathbb{C}$ (also called \textit{connection function}), the \textit{Cayley colour digraph} \cite{babai1979spectra}, denoted by $\operatorname{Cay}(G,f)$, is defined to be the directed graph with vertex set $G$ and arc set $\lbrace{(g,h)\mid g,h\in G}\rbrace$ such that each arc $(g,h)$ has colour $f(gh^{-1})$. 
The adjacency matrix of $\operatorname{Cay}(G,f)$ is defined to be the matrix whose rows and columns are indexed by the elements of $G$, and the $(g,h)$-entry is equal to $f(gh^{-1})$, i.e., $A(\operatorname{Cay}(G,f)=[f(gh^{-1})]_{g,h\in G}$. The eigenvalues of $\operatorname{Cay}(G,f)$ are the eigenvalues of its adjacency matrix $A(\operatorname{Cay}(G,f))$. Note that the adjacency matrix of $\operatorname{Cay}(G,f)$ is not symmetric, in general. In particular, when $f:G\to\lbrace{0,1}\rbrace$ and the set $S\coloneqq\lbrace{g\in G\mid f(g)=1}\rbrace$ satisfies $1\notin S$ and $S^{-1}=S$, the Cayley colour digraph $\operatorname{Cay}(G,f)$ can be identified with the Cayley graph $\operatorname{Cay}(G,S)$ and the adjacency matrix of $\operatorname{Cay}(G,f)$ coincides with that of $\operatorname{Cay}(G,S)$. 

\begin{lemma}(\cite{foster2016spectra}, Theorem $4.3$)\label{spectrum-Cay-col-digraph}
Let $f:G\to\mathbb{C}$ be a class function and $\operatorname{Irr}(G)=\lbrace{\chi_1,\ldots,\chi_m}\rbrace$. The spectrum of the Cayley colour digraph $\operatorname{Cay}(G,f)$ can be arranged as $\left\lbrace{\lambda_i^{[d_i^2]}\mid 1\le i\le m}\right\rbrace$,
where
$$\lambda_i=\frac{1}{d_i}\sum_{g\in G}f(g)\chi_i(g),
\text{ and $d_i=\chi_i(1)$ for all $1\le i\le m$}.$$
\end{lemma}

In what follows, we introduce the notion of Cayley colour graphs.

\begin{definition}
Let $G$ be a finite group and $f:G\to\mathbb{Q}$ be such that $f(g)=f(g^{-1})$ for all $g\in G$. The Cayley colour graph, denoted by $\Gamma_f=\operatorname{Cay}(G,f)$, is defined to be the undirected graph with vertex set $G$ and edge set $\lbrace{\lbrace{g,h}\rbrace\mid f(gh^{-1})\neq 0,~g,h\in G}\rbrace$, where each edge $\lbrace{g,h}\rbrace$ has some non-zero colour $f(gh^{-1})$.
\end{definition}

From the definition, it is obvious that the adjacency matrix of $\Gamma_f$ is a real symmetric matrix. Thus all the eigenvalues of $\Gamma_f$ are real. From the definition of Cayley multigraph, it is obvious that the Cayley multigraph $\operatorname{Cay}(G,S)$ is precisely the Cayley colour graph $\operatorname{Cay}(G,m_S)$. Note that if $f:G\to\mathbb{Z}_{\ge 0}$ with $f(1)=0$ and $f(g)=f(g^{-1})$ for all $g\in G$, then we can define the multiset $S_f\coloneqq\left\lbrace{g^{[f(g)]}\mid f(g)\neq 0}\right\rbrace$. In this case, the Cayley colour graph $\Gamma_f=\operatorname{Cay}(G,f)$ can be identified with the Cayley multigraph $\operatorname{Cay}(G,S_f)$. Hence for functions $f:G\to\mathbb{Z}_{\ge 0}$ with $f(1)=0$ and $f(g)=f(g^{-1})$ for all $g\in G$, we shall not distinguish between the Cayley colour graph and the Cayley multigraph. We shall denote by $\mathbb{SF}(\Gamma_f)$ and $\operatorname{Deg}(\Gamma_f)$, the splitting field and algebraic degree of the adjacency matrix of the Cayley colour graph $\Gamma_f$, respectively.

Let $C_g$ be the conjugacy class containing the element $g\in G$. Let $C_{g_1},C_{g_2},\ldots,C_{g_m}$ be the complete set of conjugacy classes of a finite group $G$. The \textit{character table} $X$ of $G$ is an $m\times m$ matrix whose $(i,j)$-entry is $\chi_i(g_j)$. By the second orthogonality relations of the irreducible characters of $G$, it can be seen that $X$ is invertible \cite{steinberg2012representation}. For a class function $f:G\to\mathbb{C}$, we define a vector relative to $f$ in $G$ by $\Delta_f=(\Delta_f(C_{g_i}))_{1\le i\le m}$, where $\Delta_f(C_{g_i})\coloneqq |C_{g_i}|f(g_i)$. This will play an important role in proving our main result in the next section. It follows immediately that, for two class functions $\alpha,\beta:G\to\mathbb{C}$ of a finite group $G$, $\Delta_\alpha=\Delta_\beta$ if and only if $\alpha=\beta$. For a function $f:G\to\mathbb{C}$ and $k\in\mathbb{Z}$, we define $f^k:G\to\mathbb{C}$ as $f^k(g)=f(g^k)$ for $g\in G$. 

\begin{lemma}\label{f-f^k-class-fn-lemma}
Let $G$ be a finite group. Then for any class function $f:G\to\mathbb{C}$ and $k\in\mathbb{Z}$, $f^k$ is also a class function of $G$.
\end{lemma}
\pf Let $g,h\in G$. Then $f^k(ghg^{-1})=f((ghg^{-1})^k)=f(gh^kg^{-1})=f(h^k)=f^k(h)$.\qed

\medskip

Let $D_n=\langle{a,b\mid a^n=b^2=1, ab=ba^{-1}}\rangle$ be the dihedral group of order $2n$. We list the character table of $D_n$ for later use.

\begin{lemma}(\cite{serre1977linear})\label{Dn-char-table-lemma}
The character table of $D_n$ is given in Table \ref{table 1}, if $n$ is odd and in Table \ref{table 2}, if $n$ is even, where $\psi_i$ and $\chi_h$ are irreducible characters of degree one and two respectively and $1\le h\le\lfloor{\frac{n-1}{2}}\rfloor$.
\end{lemma}

\begin{minipage}{0.4\linewidth}
\begin{table}[H]
\begin{center}
\begin{tabular}{lcccc}
\toprule
& & $a^k$ & & $ba^k$ \\
\midrule
$\psi_1$ & & $1$ & & $1$\\
$\psi_2$ & & $1$ & & $-1$\\
$\chi_h$ & & $2\cos{\left(\frac{2kh\pi}{n}\right)}$ & & $0$\\
\bottomrule
\end{tabular}
\caption{Character table of $D_n$ for odd $n$}
\label{table 1}
\end{center}
\end{table}
\end{minipage}
\hfill
\begin{minipage}{0.5\linewidth}
\begin{table}[H]
\begin{center}
\begin{tabular}{lcccc}
\toprule
& & $a^k$ & & $ba^k$ \\
\midrule
$\psi_1$ & & $1$ & & $1$\\
$\psi_2$ & & $1$ & & $-1$\\
$\psi_3$ & & $(-1)^k$ & & $(-1)^k$\\
$\psi_4$ & & $(-1)^k$ & & $~~~(-1)^{k+1}$\\
$\chi_h$ & & $2\cos{\left(\frac{2kh\pi}{n}\right)}$ & & $0$\\
\bottomrule
\end{tabular}
\caption{Character table of $D_n$ for even $n$}
\label{table 2}
\end{center}
\end{table}
\end{minipage}



\section{Splitting fields of Cayley colour graphs}\label{sp-field-alg-deg-Cay-col-graph-section}

From this section onward, we shall only consider Cayley colour graphs $\Gamma_f=\operatorname{Cay}(G,f)$, where $f$ is a class function. In this section, we determine the splitting field and algebraic degree of $\Gamma_f$. We first prove the following lemma.

\begin{lemma}\label{char-sum-h-Zn-lemma}
Let $G$ be a group of order $n$ and $f:G\to\mathbb{C}$ be a function. Let $\chi\in\operatorname{Irr}(G)$. Then for any $h\in\mathbb{Z}_n^*$,
$$\sum_{g\in G}f^h(g)\chi^h(g)=\sum_{g\in G}f(g)\chi(g).$$
\end{lemma}
\pf Note that every $h\in\mathbb{Z}_n^*$ permutes the elements of $G$ via the map $g\mapsto g^h$. This proves the result.\qed

\medskip

Let $G$ be a group of order $n$ and $\operatorname{Irr}(G)=\lbrace{\chi_1,\ldots,\chi_m}\rbrace$. Note that by Lemma \ref{chi-alg-int-lemma} and Lemma \ref{spectrum-Cay-col-digraph}, 
we have $\lambda_i\in\mathbb{Q}(\zeta_n)$ for all $i\in\lbrace{1,\ldots,m}\rbrace$. Thus we obtain $\mathbb{Q}\subseteq\mathbb{SF}(\Gamma_f)\subseteq\mathbb{Q}(\zeta_n)$. Suppose $K$ is a field such that $\mathbb{Q}\subseteq K\subseteq\mathbb{Q}(\zeta_n)$. Therefore, $\operatorname{Gal}(\mathbb{Q}(\zeta_n)/K)\le\operatorname{Gal}(\mathbb{Q}(\zeta_n)/\mathbb{Q})\cong\mathbb{Z}_n^{*}$. Let $\eta:\operatorname{Gal}(\mathbb{Q}(\zeta_n)/\mathbb{Q})\to\mathbb{Z}_n^{*}$ be the isomorphism such that $\sigma(\zeta_n)=\zeta_n^{\eta(\sigma)}$, where $\sigma\in\operatorname{Gal}(\mathbb{Q}(\zeta_n)/\mathbb{Q})$. Let $H_K=\eta(\operatorname{Gal}(\mathbb{Q}(\zeta_n)/K))$. Then $H_K$ is a subgroup of $\mathbb{Z}_n^{*}$. 

\begin{lemma}\label{sigma-eta-chi-lemma}
Let $G$ be a group of order $n$ and $\sigma\in\operatorname{Gal}(\mathbb{Q}(\zeta_n)/\mathbb{Q})$. Then for any character $\chi$ of $G$, $\sigma(\chi(g))=\chi(g^{\eta(\sigma)})$ for all $g\in G$.
\end{lemma}
\pf Observe that for $\sigma\in\operatorname{Gal}(\mathbb{Q}(\zeta_n)/\mathbb{Q})$, $\sigma(\gamma)=\gamma^{\eta(\sigma)}$ for any $n$th-root of unity $\gamma$. Let $\rho:G\to GL_m(\mathbb{C})$ be a representation affording the character $\chi$. Then by Lemma \ref{A-diagonalizable-lemma}, for any $g\in G$, $\rho(g)$ is diagonalizable with eigenvalues $\gamma_1,\ldots,\gamma_m$, which are $n$th-roots of unity. Thus we have,
\begin{align*}
\sigma(\chi(g))&=\sigma(\gamma_1+\cdots+\gamma_m)\\
&=\sigma(\gamma_1)+\cdots+\sigma(\gamma_m)\\
&=\gamma_1^{\eta(\sigma)}+\cdots+\gamma_m^{\eta(\sigma)}\\
&=\chi(g^{\eta(\sigma)}),
\end{align*}
where the last equality follows from the fact that $\rho(g^{\eta(\sigma)})$ is diagonalizable with eigenvalues $\gamma_1^{\eta(\sigma)},\ldots,\gamma_m^{\eta(\sigma)}$.\qed


\begin{proposition}\label{alg-integral-over-K-prop}
Let $G$ be a group of order $n$. Let $\lambda_i$ be the eigenvalues of a Cayley colour graph $\Gamma_f=\operatorname{Cay}(G,f)$ where $i\in\lbrace{1,\ldots,m}\rbrace$ and let $K$ be a field such that $\mathbb{Q}\subseteq K\subseteq\mathbb{Q}(\zeta_n)$. Then $\lambda_i\in K$ for all $i\in\lbrace{1,\ldots,m}\rbrace$ if and only if $f^h=f$ for all $h\in H_K$, where $H_K=\eta(\operatorname{Gal}(\mathbb{Q}(\zeta_n)/K))$.
\end{proposition}
\pf First we assume $f^h=f$ for all $h\in H_K$, where $H_K=\eta(\operatorname{Gal}(\mathbb{Q}(\zeta_n)/K))$. Then for any $\sigma\in\operatorname{Gal}(\mathbb{Q}(\zeta_n)/K)$, we have $\eta(\sigma)\in H_K$. Thus for any $i\in\lbrace{1,\ldots,m}\rbrace$,
\begin{align*}
\sigma(\lambda_i)&=\sigma\left(\frac{1}{d_i}\sum_{g\in G}f(g)\chi_i(g)\right)\\
&=\frac{1}{d_i}\sum_{g\in G}f(g)\sigma(\chi_i(g))\\
&=\frac{1}{d_i}\sum_{g\in G}f(g^{\eta(\sigma)})\chi_i(g^{\eta(\sigma)}) & \text{(by Lemma \ref{sigma-eta-chi-lemma})}\\
&=\frac{1}{d_i}\sum_{g\in G}f(g)\chi_i(g) & \text{(by Lemma \ref{char-sum-h-Zn-lemma})}\\
&=\lambda_i.
\end{align*}
Thus it follows that $\lambda_i\in K$ for all $i\in\lbrace{1,\ldots,m}\rbrace$.

Conversely, suppose $\lambda_i\in K$ for all $i\in\lbrace{1,\ldots,m}\rbrace$. Let $h\in H_K$ be arbitrary. Then there exists $\sigma\in\operatorname{Gal}(\mathbb{Q}(\zeta_n)/K)$ such that $\eta(\sigma)=h$. Then for all $i\in\lbrace{1,\ldots,m}\rbrace$,
\begin{align*}
\frac{1}{d_i}\sum_{g\in G}f(g)\chi_i^h(g)&=\frac{1}{d_i}\sum_{g\in G}f(g)\chi_i(g^{\eta(\sigma)})\\
&=\frac{1}{d_i}\sum_{g\in G}f(g)\sigma(\chi_i(g)) & \text{(by Lemma \ref{sigma-eta-chi-lemma})}\\
&=\sigma\left(\frac{1}{d_i}\sum_{g\in G}f(g)\chi_i(g)\right)\\
&=\sigma(\lambda_i)\\
&=\lambda_i\\
&=\frac{1}{d_i}\sum_{g\in G}f(g)\chi_i(g)\\
&=\frac{1}{d_i}\sum_{g\in G}f^h(g)\chi_i^h(g). & \text{(by Lemma \ref{char-sum-h-Zn-lemma})}
\end{align*}
Thus we have,
$$\sum_{g\in G}f(g)\chi_i^h(g)=\sum_{g\in G}f^h(g)\chi_i^h(g).$$
Again,
\begin{align*}
\sum_{g\in G}f^h(g)\chi_i^h(g)&=\sum_{g\in G}f^h(g)\chi_i(g^{\eta(\sigma)})\\
&=\sum_{g\in G}f^h(g)\sigma(\chi_i(g)) & \text{(by Lemma \ref{sigma-eta-chi-lemma})}\\
&=\sigma\left(\sum_{g\in G}f^h(g)\chi_i(g)\right),
\end{align*}
and
\begin{align*}
\sum_{g\in G}f(g)\chi_i^h(g)&=\sum_{g\in G}f(g)\chi_i(g^{\eta(\sigma)})\\
&=\sum_{g\in G}f(g)\sigma(\chi_i(g)) & \text{(by Lemma \ref{sigma-eta-chi-lemma})}\\
&=\sigma\left(\sum_{g\in G}f(g)\chi_i(g)\right).
\end{align*}
Since $\sigma$ is an automorphism, we have
$$\sum_{g\in G}f^h(g)\chi_i(g)=\sum_{g\in G}f(g)\chi_i(g).$$
Rewriting the sum over the conjugacy classes of $G$, we have
\begin{equation}\label{character-table-eigenvector}
\sum_{j=1}^m|C_{g_j}|f^h(g_j)\chi_i(g_j)=\sum_{j=1}^m|C_{g_j}|f(g_j)\chi_i(g_j).
\end{equation}
Let $X$ be the character table of $G$. By Lemma \ref{f-f^k-class-fn-lemma}, $f$ and $f^h$ are class functions of $G$. Thus from (\ref{character-table-eigenvector}) we have,
$$X\begin{bmatrix}
   \Delta_{f^h}(C_{g_1})\\
   \vdots\\
   \Delta_{f^h}(C_{g_m})
\end{bmatrix}=X\begin{bmatrix}
   \Delta_{f}(C_{g_1})\\
   \vdots\\
   \Delta_{f}(C_{g_m})
\end{bmatrix},$$
i.e.,
$$X\Delta_{f^h}=X\Delta_{f}.$$
Since $X$ is invertible, we have
$$\Delta_{f^h}=\Delta_{f.}$$
This implies $f^h=f$, which completes the proof.\qed

\begin{proposition}\label{H_f=H-prop}
Let $G$ be a group of order $n$ and let $\Gamma_f=\operatorname{Cay}(G,f)$ be a Cayley colour graph. Let $H_f\coloneqq\lbrace{h\in\mathbb{Z}_n^*\mid f^h=f}\rbrace$. Then $H_f=H$, where $H=\eta(\operatorname{Gal}(\mathbb{Q}(\zeta_n)/\mathbb{SF}(\Gamma_f)))$.
\end{proposition}
\pf Let $h\in H$, where $H=\eta(\operatorname{Gal}(\mathbb{Q}(\zeta_n)/\mathbb{SF}(\Gamma_f)))$. Since $\lambda_i\in\mathbb{SF}(\Gamma_f)$, by Proposition \ref{alg-integral-over-K-prop}, $f^h=f$. Thus $h\in H_f$.

To show the other direction, let $h'\in H_f$ and consider $\sigma=\eta^{-1}(h')$. So $h'=\eta(\sigma)$. Also we have $f^{h'}=f$. Then for any $i\in\lbrace{1,\ldots,m}\rbrace$,
\begin{align*}
\sigma(\lambda_i)&=\frac{1}{d_i}\sum_{g\in G}f(g)\sigma(\chi_i(g))\\
&=\frac{1}{d_i}\sum_{g\in G}f(g)\chi_i(g^{\eta(\sigma)}) & \text{(by Lemma \ref{sigma-eta-chi-lemma})}\\
&=\frac{1}{d_i}\sum_{g\in G}f^{h'}(g)\chi_i^{h'}(g)\\
&=\frac{1}{d_i}\sum_{g\in G}f(g)\chi_i(g) & \text{(by Lemma \ref{char-sum-h-Zn-lemma})}\\
&=\lambda_i.
\end{align*}
Thus $\sigma\in\operatorname{Gal}(\mathbb{Q}(\zeta_n)/\mathbb{SF}(\Gamma_f))$. Then $h'=\eta(\sigma)\in H$ and the result follows.\qed

\begin{theorem}\label{splitting-field-alg-deg-Cay-col-graph-thm}
Let $G$ be a group of order $n$ and let $\Gamma_f=\operatorname{Cay}(G,f)$ be a Cayley colour graph. Then the splitting field of $\Gamma_f$ is given by
$$\mathbb{SF}(\Gamma_f)=\mathbb{Q}(\zeta_n)^{\eta^{-1}(H_f)}=\lbrace{x\in\mathbb{Q}(\zeta_n)\mid\sigma(x)=x,\text{for all } \sigma\in\eta^{-1}(H_f)}\rbrace.$$
Moreover, the algebraic degree of $\Gamma_f$ is
$$\operatorname{Deg}(\Gamma_f)=\frac{\varphi(n)}{|H_f|},$$
where $\varphi(\cdot)$ is the Euler totient function and $H_f=\lbrace{h\in\mathbb{Z}_n^*\mid f^h=f}\rbrace$.
\end{theorem}
\pf Since $\mathbb{Q}(\zeta_n)/\mathbb{SF}(\Gamma_f)$ is a Galois extension, by Proposition \ref{H_f=H-prop} we get $\mathbb{SF}(\Gamma_f)=\mathbb{Q}(\zeta_n)^{\eta^{-1}(H_f)}$. It follows that 
$$\operatorname{Deg}(\Gamma_f)=[\mathbb{SF}(\Gamma_f):\mathbb{Q}]=\frac{[\mathbb{Q}(\zeta_n):\mathbb{Q}]}{[\mathbb{Q}(\zeta_n):\mathbb{SF}(\Gamma_f)]}=\frac{\varphi(n)}{|H_f|}.$$\qed

\medskip

By Theorem \ref{splitting-field-alg-deg-Cay-col-graph-thm}, we have the following characterization of rational Cayley colour graphs.

\begin{corollary}\label{Cay-col-graph-rational-iff-cond-corollary}
Let $G$ be a group of order $n$. Then the Cayley colour graph $\Gamma_f=\operatorname{Cay}(G,f)$ is rational if and only if $f^h=f$ for all $h\in\mathbb{Z}_n^{*}$. 
\end{corollary}

\begin{corollary}\label{Cay-col-graph-integral-iff-cond-corollary}
If $f:G\to\mathbb{Q}$ be such that $f(G)\subseteq\mathbb{Z}$, then the Cayley colour graph $\Gamma_f=\operatorname{Cay}(G,f)$ is integral if and only if $f^h=f$ for all $h\in\mathbb{Z}_n^{*}$. 
\end{corollary}

\begin{corollary}\label{alg-deg-Cay-col-graph-divides-corollary}
Let $G$ be a group of order $n\ge 3$ and let $\Gamma_f=\operatorname{Cay}(G,f)$ be a Cayley colour graph. Then $\operatorname{Deg}(\Gamma_f)$ divides $\frac{\varphi(n)}{2}$.
\end{corollary}
\pf Clearly, $f^1=f$. Also, since $f(g^{-1})=f(g)$ for all $g\in G$, we have $f^{-1}=f$. Hence $\lbrace{\pm 1}\rbrace\le H_f$, i.e., $H_f$ must be of even order.\qed 

\medskip

We provide two examples to illustrate Theorem \ref{splitting-field-alg-deg-Cay-col-graph-thm} and Corollary \ref{Cay-col-graph-integral-iff-cond-corollary}.

\begin{example}
Let $G=D_8$. Let $\alpha:G\to\mathbb{Q}$ be defined as
$$\alpha(1)=\alpha(a^4)=0,$$ $$\alpha(a)=\alpha(a^7)=1,$$
$$\alpha(a^2)=\alpha(a^6)=1/2,$$
$$\alpha(a^3)=\alpha(a^5)=3/5,$$
$$\alpha(b)=\alpha(ba^2)=\alpha(ba^4)=\alpha(ba^6)=4,$$
$$\alpha(ba)=\alpha(ba^3)=\alpha(ba^5)=\alpha(ba^7)=7.$$
Then $\alpha$ is a class function satisfying $\alpha(g)=\alpha(g^{-1})$ for all $g\in G$. Let $\Gamma_\alpha=\operatorname{Cay}(G,\alpha)$. Then it follows that $H_\alpha=\lbrace{\pm1,\pm7}\rbrace\le\mathbb{Z}_{16}^{*}$. By Theorem \ref{splitting-field-alg-deg-Cay-col-graph-thm}, we have $\mathbb{SF}(\Gamma_\alpha)=\mathbb{Q}(\zeta_{16})^{\eta^{-1}(H_\alpha)}=\mathbb{Q}(\sqrt{2})$ and $\operatorname{Deg}(\Gamma_\alpha)=\frac{\varphi(16)}{|H_\alpha|}=2$. Indeed by Lemma \ref{spectrum-Cay-col-digraph} and \ref{Dn-char-table-lemma}, we have
$$\operatorname{Spec}(\Gamma_\alpha)=\begin{pmatrix}
\frac{241}{5} & \frac{49}{5} & -\frac{71}{5} & -\frac{199}{5} & -1 & \frac{2\sqrt{2}}{5} & -\frac{2\sqrt{2}}{5}\\
1 & 1 & 1 & 1 & 4 & 4 & 4
\end{pmatrix}.$$
\end{example}

\begin{example}
Let $G=D_8$. Let $\beta:G\to\mathbb{Q}$ be defined as
$$\beta(1)=\beta(a^4)=0,$$
$$\beta(a)=\beta(a^3)=\beta(a^5)=\beta(a^7)=1,$$
$$\beta(a^2)=\beta(a^6)=3,$$
$$\beta(b)=\beta(ba^2)=\beta(ba^4)=\beta(ba^6)=4,$$
$$\beta(ba)=\beta(ba^3)=\beta(ba^5)=\beta(ba^7)=8.$$
Then $\beta$ is a class function satisfying $\beta(g)=\beta(g^{-1})$ for all $g\in G$. Let $\Gamma_\beta=\operatorname{Cay}(G,\beta)$. It can be seen that $\beta^h=\beta$ for all $h\in\mathbb{Z}_{16}^{*}$. Hence by Corollary \ref{Cay-col-graph-integral-iff-cond-corollary}, $\Gamma_\beta$ is integral. Indeed by Lemma \ref{spectrum-Cay-col-digraph} and \ref{Dn-char-table-lemma}, we have 
$$\operatorname{Spec}(\Gamma_\beta)=\begin{pmatrix}
58 & 18 & -14 & -38 & -6 & 0\\
1 & 1 & 1 & 1 & 4 & 8
\end{pmatrix}.$$

\end{example}

\section{Algebraic degree and distance algebraic degree of normal Cayley graphs}\label{alg-deg-dist-alg-deg-normal-Cay-graph-section}

We now consider our conventional Cayley (multi)graphs $\operatorname{Cay}(G,S)$. In \cite{zhang20262}, Zhang et al. have determined the splitting field and algebraic degree, along with the distance splitting field and distance algebraic degree of a normal Cayley graph. In this section, we do the same but as a corollary of Theorem \ref{splitting-field-alg-deg-Cay-col-graph-thm}. For any $S\subseteq G$ and $\chi\in\operatorname{Irr}(G)$, define
$$\chi(S)\coloneqq\sum_{s\in S}\chi(s).$$
For any (multi)set $S\subseteq G$ and $k\in\mathbb{Z}$, define 
$$S^k\coloneqq\lbrace{s^k\mid s\in S}\rbrace.$$

\subsection{\textbf{Algebraic degree of normal Cayley multigraphs}}

\medskip

Let $\Gamma=\operatorname{Cay}(G,S)$ be a Cayley multigraph. Recall the multiplicity function $m_S:G\to\mathbb{Z}_{\ge 0}$ defined as
$$m_S(x)=\begin{cases}
    \#\text{ $x$ that occurs in S}, & \text{if $x\in S$},\\
    0, & \text{if $x\notin S$.}
\end{cases}$$
Since $S^{-1}=S$, $m_S(g)=m_S(g^{-1})$ for all $g\in G$. If $S$ is normal, it follows that $m_S$ is a class function. Thus the adjacency matrix of $\Gamma$ is identified with that of the Cayley colour graph $\Gamma_{m_S}=\operatorname{Cay}(G,m_S)$. From Lemma \ref{spectrum-Cay-col-digraph}, the spectrum of $\Gamma$ can be obtained as the following.

\begin{lemma}\label{normal-Cay-multigraph-spectrum-lemma}
Let $G$ be a finite group and $\operatorname{Irr}(G)=\lbrace{\chi_1,\ldots,\chi_m}\rbrace$. Let $S\subseteq G\setminus\lbrace{1}\rbrace$ be a normal connection multiset. Then the spectrum of the normal Cayley multigraph $\Gamma=\operatorname{Cay}(G,S)$ be arranged as $\left\lbrace{\theta_i^{[d_i^2]}\mid 1\le i\le m}\right\rbrace$,
where
$$\theta_i=\frac{1}{d_i}\sum_{x\in S}m_S(x)\chi_i(x), \text{ and $d_i=\chi_i(1)$ for all $1\le i\le m$}.$$
\end{lemma}

Recall from Proposition \ref{H_f=H-prop} that for a class function $f$, $H_f=\lbrace{h\in\mathbb{Z}_n^{*}\mid f^h=f}\rbrace$. By $S^h=S$, we shall mean that the equality holds as multisets.

\begin{lemma}\label{H-del-S-H*-lemma}
Let $H^*\coloneqq\lbrace{h\in\mathbb{Z}_n^{*}\mid S^h=S}\rbrace$. Then $H_{m_S}=H^*$.
\end{lemma}
\pf Let $h\in H_{m_S}$. Then $m_S^h=m_S$. Let $x\in S$. Then $m_S(x)=m_S^h(x)=m_S(x^h)$. Now let $x\notin S$. Then $m_S(x)=m_S^h(x)=m_S(x^h)=0$, i.e., $x^h\notin S$. Hence for any $x\in G$, $x\in S$ if and only if $x^h\in S$ with $m_S(x)=m_S(x^h)$, which implies $S^h=S$ as multisets. Thus $h\in H^*$.

To show the other direction, let $h^*\in H^*$. Then $S^{h^*}=S$ as multisets. Let $x\in G$. Then $x\in S$ if and only if $x^{h^*}\in S$ with $m_S(x)=m_S(x^{h^*})$. Hence $m_S^{h^*}=m_S$, which implies that $h^*\in H_{m_S}$.\qed

\medskip

From Lemma \ref{H-del-S-H*-lemma} and Theorem \ref{splitting-field-alg-deg-Cay-col-graph-thm}, we can determine the splitting field and algebraic degree of the normal Cayley multigraph $\Gamma$.

\begin{corollary}\label{sp-field-alg-deg-normal-Cay-multigraph-corollary}
Let $G$ be a group of order $n$ and $\Gamma=\operatorname{Cay}(G,S)$ be a normal Cayley multigraph. Then the splitting field of $\Gamma$ is
$$\mathbb{SF}(\Gamma)=\mathbb{Q}(\zeta_n)^{\eta^{-1}(H^*)}=\lbrace{x\in\mathbb{Q}(\zeta_n)\mid\sigma(x)=x,\text{for all } \sigma\in\eta^{-1}(H^*)}\rbrace.$$
Moreover, the algebraic degree of $\Gamma$ is
$$\operatorname{Deg}(\Gamma)=\frac{\varphi(n)}{|H^*|},$$
where $H^*=\lbrace{h\in\mathbb{Z}_n^*\mid S^h=S}\rbrace$.
\end{corollary}

\begin{corollary}\label{normal-Cay-multigraph-intgeral-corollary}
Let $G$ be a group of order $n$. Then the normal Cayley multigraph $\Gamma=\operatorname{Cay}(G,S)$ is integral if and only if $S^h=S$ as multisets, for all $h\in\mathbb{Z}_n^{*}$. 
\end{corollary}

If we consider $\Gamma=\operatorname{Cay}(G,S)$ as a normal Cayley graph, then we have the exact same expressions for the splitting field and algebraic degree of $\Gamma$ [\cite{wu2024splitting}, Theorem $3.1$]. In that case, $S^h=S$ will be considered as set equality. Thus we have

\begin{corollary}(\cite{wu2024splitting}, Corollary $3.1$)
Let $G$ be a group of order $n$. Then the normal Cayley graph $\Gamma=\operatorname{Cay}(G,S)$ is integral if and only if $S^h=S$ for all $h\in\mathbb{Z}_n^{*}$. 
\end{corollary}

It follows immediately that for $h\in\mathbb{Z}_n^{*}$, if $S^h=S$ holds as multisets, then $S^h=S$ holds as sets also. This leads us to the following result.

\begin{corollary}
Let $\Gamma=\operatorname{Cay}(G,S)$ be a normal Cayley multigraph and $\overline{\Gamma}=\operatorname{Cay}(G,\overline{S})$ be the normal Cayley graph, where $\overline{S}$ is obtained from $S$ by removing the repeated elements. Then $\operatorname{Deg}(\overline{\Gamma})$ divides $\operatorname{Deg}(\Gamma)$.
\end{corollary}

We provide two examples to illustrate Corollary \ref{sp-field-alg-deg-normal-Cay-multigraph-corollary} and \ref{normal-Cay-multigraph-intgeral-corollary}.

\begin{example}
Let $G=D_5$ and consider the multiset $S_1=\lbrace{a^2,a^2,a^3,a^3,b,ba,ba^2,ba^3,ba^4}\rbrace$. Then $\Gamma_1=\operatorname{Cay}(G,S_1)$ is a normal Cayley multigraph. It follows that $H^{*}=\lbrace{\pm1}\rbrace\le\mathbb{Z}_{10}^{*}$. By Corollary \ref{sp-field-alg-deg-normal-Cay-multigraph-corollary}, we have $\mathbb{SF}(\Gamma_1)=\mathbb{Q}(\zeta_{10})^{\eta^{-1}(H^{*})}=\mathbb{Q}(\sqrt{5})$ and $\operatorname{Deg}(\Gamma_1)=\frac{\varphi(10)}{|H^{*}|}=2$. Indeed by Lemma \ref{normal-Cay-multigraph-spectrum-lemma}  and \ref{Dn-char-table-lemma}, we have
$$\operatorname{Spec}(\Gamma_1)=\begin{pmatrix}
9 & -1 & -1+\sqrt{5} & -1-\sqrt{5}\\
1 & 1 & 4 & 4
\end{pmatrix}.$$
\end{example}

\begin{example}
Let $G=D_5$ and consider the normal multiset $S_2=\lbrace{a,a,a^2,a^2,a^3,a^3,a^4,a^4}\rbrace$. Let $\Gamma_2=\operatorname{Cay}(G,S_2)$. Then it is easy to see that $S_2^h=S_2$ as multisets for all $h\in\mathbb{Z}_{10}^{*}$. Thus by Corollary \ref{normal-Cay-multigraph-intgeral-corollary}, $\Gamma_2$ is integral. Indeed by Lemma \ref{normal-Cay-multigraph-spectrum-lemma} and \ref{Dn-char-table-lemma}, we have
$$\operatorname{Spec}(\Gamma_2)=\begin{pmatrix}
8 & -2\\
2 & 8
\end{pmatrix}.$$
\end{example}


\subsection{\textbf{Distance algebraic degree of normal Cayley graphs}}

\medskip

Let $\Gamma=\operatorname{Cay}(G,S)$ be a connected Cayley graph. Consider the function $\ell_S:G\to\mathbb{C}$ defined as
$$\ell_S(x)=\begin{cases}
    \operatorname{min}\lbrace{k\mid x=s_1\cdots s_k,\text{ for $s_i\in S$}}\rbrace, & \text{ if $x\neq 1$},\\
    0, & \text{ if $x=1$}.
\end{cases}$$

In terms of distance, the function $\ell_S$ satisfies $\ell_S(g)=d_\Gamma(g,1)$ and hence $\ell_S(g)=\ell_S(g^{-1})$ for all $g\in G$. It can be seen from [\cite{wu2024splitting}, Lemma $4.1$] that if $S$ is normal, then $\ell_S$ is a class function. Thus the distance matrix of $\Gamma$ is identified with the adjacency matrix of the Cayley colour graph $\Gamma_{\ell_S}=\operatorname{Cay}(G,\ell_S)$. From Lemma \ref{spectrum-Cay-col-digraph}, the distance spectrum of $\Gamma$ can be obtained as the following.

\begin{lemma}
Let $G$ be a finite group and $\operatorname{Irr}(G)=\lbrace{\chi_1,\ldots,\chi_m}\rbrace$. Let $S\subseteq G\setminus\lbrace{1}\rbrace$ be a normal connection set. Then the distance spectrum of the normal Cayley graph $\Gamma=\operatorname{Cay}(G,S)$ be arranged as $\left\lbrace{\mu_i^{[d_i^2]}\mid 1\le i\le m}\right\rbrace$,
where
$$\mu_i=\frac{1}{d_i}\sum_{g\in G}\ell_S(g)\chi_i(g), \text{ and $d_i=\chi_i(1)$ for all $1\le i\le m$}.$$
\end{lemma}

Let $\Gamma=\operatorname{Cay}(G,S)$ be a connected, normal Cayley graph. The $i$-th layer of $\Gamma$ with respect to the vertex $1$, denoted
by $S_i$, is the set of vertices of  whose distance to the vertex $1$ is exactly $i$. Clearly, $S_0=\lbrace{1}\rbrace$ and $S_1=S$. Let $d$ be the diameter of $\Gamma$. Then $S_0,S_1,\ldots,S_d$ is a partition of $G$ and $\ell_S(g)=i$ if $g\in S_i$. Thus the distance eigenvalue $\mu_i$ can be written in the following form:

$$\mu_i=\frac{1}{d_i}\left(\chi_i(S_1)+2\chi_i(S_2)+\cdots+d\chi_i(S_d)\right).$$

\begin{lemma}\label{H-ell-S-H*-lemma}
Let $H'\coloneqq\lbrace{h\in\mathbb{Z}_n^{*}\mid S_i^h=S_i,~1\le i\le d}\rbrace$. Then $H_{\ell_S}=H'$.
\end{lemma}
\pf Let $h\in H_{\ell_S}$. Then $\ell_S^h=\ell_S$. Let $i\in\lbrace{1,\ldots,d}\rbrace$ be fixed and let $x\in S_i$. Then $\ell_S(x)=\ell_S^h(x)=\ell_S(x^h)=i$, which implies $x^h\in S_i$. Thus $S_i^h\subseteq S_i$. As $|S_i^h|=|S_i|$, we have $S_i^h=S_i$, for all $i\in\lbrace{1,\ldots,d}\rbrace$. Hence $h\in H'$.

To show the other direction, let $h'\in H'$. Then $S_i^{h'}=S_i$ for all $i\in\lbrace{1,\ldots,d}\rbrace$. Let $x\in G$. Since $\lbrace{S_i\mid 1\le i\le d}\rbrace$ partitions $G$, there exists a unique $j\in\lbrace{1,\ldots,d}\rbrace$ such that $x\in S_j$. Then $\ell_S^{h'}(x)=\ell_S(x^{h'})=\ell_S(x)=i$. Thus $\ell_S^{h'}=\ell_S$, which implies $h'\in H_{\ell_S}$.\qed

\medskip

From Lemma \ref{H-ell-S-H*-lemma} and Theorem \ref{splitting-field-alg-deg-Cay-col-graph-thm}, we can determine the distance splitting field and distance algebraic degree of the normal Cayley graph $\Gamma$.

\begin{corollary}\label{dist-sp-field-dist-alg-deg-normal-Cay-graph-corollary}(\cite{wu2024splitting}, Theorem $4.1$)
Let $G$ be a group of order $n$ and $\Gamma=\operatorname{Cay}(G,S)$ be a normal Cayley graph with diameter $d$. Let $S_i$ be the set of vertices of $\Gamma$ whose distance to the vertex $1$ is exactly $i$. Then the distance splitting field of $\Gamma$ is
$$\mathbb{SF}_D(\Gamma)=\mathbb{Q}(\zeta_n)^{\eta^{-1}(H')}=\lbrace{x\in\mathbb{Q}(\zeta_n)\mid\sigma(x)=x,\text{for all } \sigma\in\eta^{-1}(H')}\rbrace.$$
Moreover, the distance algebraic degree of $\Gamma$ is
$$\operatorname{Deg}_D(\Gamma)=\frac{\varphi(n)}{|H'|},$$
where $H'=\lbrace{h\in\mathbb{Z}_n^*\mid S_i^h=S_i,~1\le i\le d}\rbrace$.
\end{corollary}

\begin{corollary}(\cite{wu2024splitting}, Corollary $4.1$)
Let $G$ be a group of order $n$. Then the normal Cayley graph $\Gamma=\operatorname{Cay}(G,S)$ is distance integral if and only if $S_i^h=S_i$ for all $1\le i\le d$ and $h\in\mathbb{Z}_n^{*}$.
\end{corollary}

From Corollary \ref{alg-deg-Cay-col-graph-divides-corollary},\ref{sp-field-alg-deg-normal-Cay-multigraph-corollary} and \ref{dist-sp-field-dist-alg-deg-normal-Cay-graph-corollary}, we can deduce the following simple fact.

\begin{corollary}
Let $G$ be a group of order $n$ and let $\Gamma=\operatorname{Cay}(G,S)$ be a normal Cayley graph. Then $\operatorname{Deg}(\Gamma)$ and $\operatorname{Deg}_D(\Gamma)$ divides $\frac{\varphi(n)}{2}$.
\end{corollary}


\section{Algebraic integrality of Cayley colour graphs over a field $K$}\label{alg-int-Cay-col-graph-section}

Recall that for a field $K$ with $\mathbb{Q}\subseteq K\subseteq\mathbb{Q}(\zeta_n)$, we call a Cayley colour graph $\Gamma_f$ to be algebraically integral over $K$, if all the eigenvalues of $\Gamma_f$ are in $K$. In this section, we find the relation between the algebraic integrality of two Cayley colour graphs. Let $\Gamma_\alpha=\operatorname{Cay}(G,\alpha)$ and $\Gamma_\beta=\operatorname{Cay}(G,\beta)$ be two Cayley colour graphs. Recall from Proposition \ref{alg-integral-over-K-prop} that $H_K=\eta(\operatorname{Gal}(\mathbb{Q}(\zeta_n)/K))$. Define,
$$H_\alpha(K)=\lbrace{h\in H_K\mid \alpha^h=\alpha}\rbrace,$$
$$H_\beta(K)=\lbrace{h\in H_K\mid \beta^h=\beta}\rbrace.$$

\begin{theorem}\label{algebraically-integral-over-K-thm}
Let $H_\alpha(K)=H_\beta(K)$. Then $\Gamma_\alpha$ is algebraically integral over $K$ if and only if $\Gamma_\beta$ is algebraically integral over $K$.

\end{theorem}
\pf 
We only prove the necessity. The sufficiency follows in the exact manner. Suppose $\Gamma_\alpha$ is algebraically integral over $K$. Then all the eigenvalues of $\Gamma_\alpha$ are in $K$. By Proposition \ref{alg-integral-over-K-prop}, $\alpha^h=\alpha$ for all $h\in H_K$, where $H_K=\eta(\operatorname{Gal}(\mathbb{Q}(\zeta_n)/K))$. This implies $H_\alpha(K)=H_K$. By the hypothesis we have, $H_\beta(K)=H_K$. Hence $\beta^h=\beta$ for all $h\in H_K$, where $H_K=\eta(\operatorname{Gal}(\mathbb{Q}(\zeta_n)/K))$. By Proposition \ref{alg-integral-over-K-prop}, all the eigenvalues of $\Gamma_\beta$ are in $K$. Hence $\Gamma_\beta$ is algebraically integral over $K$.\qed

\begin{corollary}\label{same-alg-deg-Cay-col-graph-corollary}
Let $H_\alpha(K)=H_\beta(K)$, for any field $K$ with $\mathbb{Q}\subseteq K\subseteq\mathbb{Q}(\zeta_n)$. Then $\operatorname{Deg}(\Gamma_\alpha)=\operatorname{Deg}(\Gamma_\beta)$.
\end{corollary}
\pf It is enough to show that the splitting fields of $\Gamma_\alpha$ and $\Gamma_\beta$ are same. We already have $\mathbb{Q}\subseteq\mathbb{SF}(\Gamma_\alpha),\mathbb{SF}(\Gamma_\beta)\subseteq\mathbb{Q}(\zeta_n)$. Since $\Gamma_\alpha$ is algebraically integral over $\mathbb{SF}(\Gamma_\alpha)$, by Theorem \ref{algebraically-integral-over-K-thm}, $\Gamma_\beta$ is algebraically integral over $\mathbb{SF}(\Gamma_\alpha)$. Hence all the eigenvalues of $\Gamma_\beta$ are in $\mathbb{SF}(\Gamma_\alpha)$. Thus $\mathbb{SF}(\Gamma_\beta)\subseteq\mathbb{SF}(\Gamma_\alpha)$. 

Conversely, Since $\Gamma_\beta$ is algebraically integral over $\mathbb{SF}(\Gamma_\beta)$, by Theorem \ref{algebraically-integral-over-K-thm}, $\Gamma_\alpha$ is algebraically integral over $\mathbb{SF}(\Gamma_\beta)$. Hence all the eigenvalues of $\Gamma_\alpha$ are in $\mathbb{SF}(\Gamma_\beta)$. Thus $\mathbb{SF}(\Gamma_\alpha)\subseteq\mathbb{SF}(\Gamma_\beta)$.\qed

\medskip

In [\cite{zhang20262}, Proposition $1$], the authors have shown that for any $\mathbb{Q}\subseteq K\subseteq\mathbb{Q}(\zeta_n)$, the normal Cayley graph $\Gamma=\operatorname{Cay}(G,S)$ is algebraically integral over $K$ if and only if it is distance algebraically integral over $K$. Following the same argument as in Corollary \ref{same-alg-deg-Cay-col-graph-corollary}, we can show that $\mathbb{SF}(\Gamma)=\mathbb{SF}_D(\Gamma)$ and hence $\operatorname{Deg}(\Gamma)=\operatorname{Deg}_D(\Gamma)$. Thus we have the following result, generalizing [\cite{zhang20262}, Theorem $1(1)$], which states that for any connected graph $\Gamma=\operatorname{Cay}(G,S)$ on a cyclic group $G$, $\Gamma$ is $2$-integral if and only if $\Gamma$ is $2$-distance integral.

\begin{corollary}\label{normal-Cay-d-int-iff-d-dist-int-corollary}
Let $\Gamma=\operatorname{Cay}(G,S)$ be a normal Cayley graph. Then $\Gamma$ is $d$-integral if and only if $\Gamma$ is $d$-distance integral.
\end{corollary}

\begin{remark}
If we drop the normality condition on $\Gamma$, then the result does not hold in general, as shown by Zhang et al. [\cite{zhang20262}, Example $2$ and Example $3$].
\end{remark}

\begin{remark}
It was proved in [\cite{poddar2025non}, Proposition $2.3$], that for every positive divisor $d$ of $\frac{\varphi(n)}{2}$, $(n\ge 3)$, there exists a $\frac{\varphi(n)}{d}$-regular circulant graph of order $n$ with algebraic degree $d$. From Corollary \ref{normal-Cay-d-int-iff-d-dist-int-corollary}, we get the same circulant graph with distance algebraic degree $d$. However, the same is not true for arbitrary normal Cayley graph $\Gamma$. For example, by exhaustive search, one may verify that although $2$ divides $\frac{\varphi(8)}{2}$, there does not exist a $2$-integral (and hence $2$-distance integral) normal Cayley graph on $D_4$.
\end{remark}

\section*{Acknowledgement}
The author is supported by the funding of UGC [NTA Ref. No. 211610129182], Govt. of India. The author is grateful to his supervisor Angsuman Das and also to Dr. Sucharita Biswas, Indian Institute of Technology Bombay, India for many fruitful discussions.

\bibliographystyle{abbrv}
\bibliography{ref}

@book{serre1977linear,
  title={Linear representations of finite groups},
  author={Serre, Jean Pierre},
  volume={42},
  year={1977},
  publisher={Springer}
}

@book{steinberg2012representation,
  title={Representation theory of finite groups: an introductory approach},
  author={Steinberg, Benjamin},
  year={2012},
  publisher={Springer}
}

@book{isaacs1994character,
  title={Character theory of finite groups},
  author={Isaacs, I Martin},
  volume={69},
  year={1994},
  publisher={Courier Corporation}
}

@inproceedings{harary1974graphs,
  title={Which graphs have integral spectra?},
  author={Harary, Frank and Schwenk, Allen J},
  booktitle={Graphs and Combinatorics: Proceedings of the Capital Conference on Graph Theory and Combinatorics at the George Washington University June 18--22, 1973},
  pages={45--51},
  year={1974},
  organization={Springer}
}

@article{so2006integral,
  title={Integral circulant graphs},
  author={So, Wasin},
  journal={Discrete Mathematics},
  volume={306},
  number={1},
  pages={153--158},
  year={2006},
  publisher={Elsevier}
}

@article{babai1979spectra,
  title={Spectra of Cayley graphs},
  author={Babai, L{\'a}szl{\'o}},
  journal={Journal of Combinatorial Theory, Series B},
  volume={27},
  number={2},
  pages={180--189},
  year={1979},
  publisher={Elsevier}
}

@article{monius2018graphs,
  title={Which graphs have non-integral spectra?},
  author={M{\"o}nius, Katja and Steuding, J{\"o}rn and Stumpf, Pascal},
  journal={Graphs and Combinatorics},
  volume={34},
  pages={1507--1518},
  year={2018},
  publisher={Springer}
}

@article{monius2020algebraic,
  title={The algebraic degree of spectra of circulant graphs},
  author={M{\"o}nius, Katja},
  journal={Journal of Number Theory},
  volume={208},
  pages={295--304},
  year={2020},
  publisher={Elsevier}
}

@article{monius2022splitting,
  title={Splitting fields of spectra of circulant graphs},
  author={M{\"o}nius, Katja},
  journal={Journal of Algebra},
  volume={594},
  pages={154--169},
  year={2022},
  publisher={Elsevier}
}

@article{lu2023algebraic,
  title={Algebraic degree of Cayley graphs over abelian groups and dihedral groups},
  author={Lu, Lu and M{\"o}nius, Katja},
  journal={Journal of Algebraic Combinatorics},
  volume={57},
  number={3},
  pages={753--761},
  year={2023},
  publisher={Springer}
}

@article{abdollahi20242,
  title={On $2 $-integral Cayley graphs},
  author={Abdollahi, Alireza and Arezoomand, Majid and Feng, Tao and Wang, Shixin},
  journal={arXiv preprint arXiv:2401.15306},
  year={2024}
}

@article{bridges1982rational,
  title={Rational G-matrices with rational eigenvalues},
  author={Bridges, William G and Mena, Roberto A},
  journal={Journal of Combinatorial Theory, Series A},
  volume={32},
  number={2},
  pages={264--280},
  year={1982},
  publisher={Elsevier}
}

@article{alperin2012integral,
  title={Integral sets and Cayley graphs of finite groups},
  author={Alperin, Roger C and Peterson, Brian L},
  journal={The Electronic Journal of Combinatorics},
  pages={1--12},
  year={2012}
}

@article{huang2021integral,
  title={Integral and distance integral Cayley graphs over generalized dihedral groups},
  author={Huang, Jing and Li, Shuchao},
  journal={Journal of Algebraic Combinatorics},
  volume={53},
  number={4},
  pages={921--943},
  year={2021},
  publisher={Springer}
}

@article{cheng2019integral,
  title={Integral Cayley graphs over dicyclic group},
  author={Cheng, Tao and Feng, Lihua and Huang, Hualin},
  journal={Linear Algebra and its Applications},
  volume={566},
  pages={121--137},
  year={2019},
  publisher={Elsevier}
}

@article{lu2018integral,
  title={Integral Cayley graphs over dihedral groups},
  author={Lu, Lu and Huang, Qiongxiang and Huang, Xueyi},
  journal={Journal of algebraic combinatorics},
  volume={47},
  pages={585--601},
  year={2018},
  publisher={Springer}
}

@article{cheng2023integral,
  title={Integral Cayley graphs over semi-dihedral groups},
  author={Cheng, Tao and Feng, Lihua and Yu, Guihai and Zhang, Chi},
  journal={Applicable Analysis and Discrete Mathematics},
  volume={17},
  number={2},
  pages={334--356},
  year={2023},
  publisher={JSTOR}
}

@article{godsil2025integral,
  title={Integral normal Cayley graphs},
  author={Godsil, Chris and Spiga, Pablo},
  journal={Journal of Algebraic Combinatorics},
  volume={62},
  number={1},
  pages={20},
  year={2025},
  publisher={Springer}
}

@article{guo2019integral,
  title={Integral Cayley graphs},
  author={Guo, W and Lytkina, Daria Viktorovna and Mazurov, Victor Danilovich and Revin, Danila Olegovich},
  journal={Algebra and Logic},
  volume={58},
  number={4},
  pages={297--305},
  year={2019},
  publisher={Springer}
}

@article{poddar2025non,
  title={Non-isomorphic $ d $-integral circulant graphs},
  author={Poddar, Sauvik and Das, Angsuman},
  journal={arXiv preprint arXiv:2507.17407},
  year={2025}
}

@article{wu2024splitting,
  title={Splitting fields of some matrices of normal (mixed) Cayley graphs},
  author={Wu, Yongjiang and Guo, Qinghong and Yang, Jing and Feng, Lihua},
  journal={Discrete Mathematics},
  volume={347},
  number={5},
  pages={113914},
  year={2024},
  publisher={Elsevier}
}

@article{zhang20262,
  title={On 2-distance integral Cayley graphs},
  author={Zhang, Xiaoqian and Liu, Weijun and Lu, Lu},
  journal={Applied Mathematics and Computation},
  volume={517},
  pages={129893},
  year={2026},
  publisher={Elsevier}
}

@article{liu2022eigenvalues,
  title={Eigenvalues of Cayley Graphs},
  author={Liu, Xiaogang and Zhou, Sanming},
  journal={The Electronic Journal of Combinatorics},
  pages={P2--9},
  year={2022}
}

@article{li2013circulant,
  title={Circulant digraphs integral over number fields},
  author={Li, Fei},
  journal={Discrete Mathematics},
  volume={313},
  number={6},
  pages={821--823},
  year={2013},
  publisher={Elsevier}
}

@article{li2020method,
  title={A method to determine algebraically integral Cayley digraphs on finite abelian group},
  author={Li, Fei},
  journal={Contributions to Discrete Mathematics},
  volume={15},
  number={2},
  pages={148--152},
  year={2020}
}

@article{huang2021distance,
  title={Distance-integral Cayley graphs over abelian groups and dicyclic groups},
  author={Huang, Jing and Li, Shuchao},
  journal={Journal of Algebraic Combinatorics},
  volume={54},
  number={4},
  pages={1047--1063},
  year={2021},
  publisher={Springer}
}

@article{foster2016spectra,
  title={Spectra of Cayley graphs of complex reflection groups},
  author={Foster-Greenwood, Briana and Kriloff, Cathy},
  journal={Journal of Algebraic Combinatorics},
  volume={44},
  number={1},
  pages={33--57},
  year={2016},
  publisher={Springer}
}

@article{devos2013integral,
  title={Integral Cayley Multigraphs over Abelian and Hamiltonian Groups},
  author={DeVos, Matt and Krakovski, Roi and Mohar, Bojan and Ahmady, Azhvan Sheikh},
  journal={The Electronic Journal of Combinatorics},
  pages={P63--P63},
  year={2013}
}

@article{sripaisan2022algebraic,
  title={Algebraic degree of spectra of Cayley hypergraphs},
  author={Sripaisan, Naparat and Meemark, Yotsanan},
  journal={Discrete Applied Mathematics},
  volume={316},
  pages={87--94},
  year={2022},
  publisher={Elsevier}
}

@article{wang2024algebraic,
  title={Algebraic degrees of quasi-abelian semi-Cayley digraphs},
  author={Wang, Shixin and Arezoomand, Majid and Feng, Tao},
  journal={Discrete Mathematics},
  volume={347},
  number={12},
  pages={114178},
  year={2024},
  publisher={Elsevier}
}

@article{li2025algebraic,
  title={Algebraic degrees of n-Cayley digraphs over abelian groups},
  author={Li, Hao and Liu, Xiaogang},
  journal={Journal of Algebraic Combinatorics},
  volume={61},
  number={3},
  pages={40},
  year={2025},
  publisher={Springer}
}

\end{document}